\documentclass[a4paper, 11pt]{article}

\usepackage[font={scriptsize}]{caption}
\usepackage{subcaption}

\usepackage{nomencl}
\makenomenclature

\usepackage{etoolbox}
\patchcmd{\thebibliography}{\section*{\refname}}{}{}{}

\let\oldbibliography\thebibliography
\renewcommand{\thebibliography}[1]{\oldbibliography{#1}
	\setlength{\itemsep}{-4pt}} 

\usepackage{scrpage2}
\ifoot[]{}
\ofoot[\pagemark]{\pagemark}

\usepackage{graphicx,amsmath,amssymb}
\usepackage{epstopdf}

\usepackage{etoolbox}
\usepackage{lmodern}
\usepackage[affil-it]{authblk} 

\usepackage[scaled]{helvet}
\renewcommand\familydefault{\sfdefault} 
\usepackage[T1]{fontenc}

\newcommand{\lmr}{\fontfamily{lmr}\selectfont} 

\usepackage{titling}
\usepackage{tgbonum}

\usepackage{paralist}
\usepackage{color}

\usepackage{geometry}
\usepackage{url}

\newcommand\blfootnote[1]{%
	\begingroup
	\renewcommand\thefootnote{}\footnote{#1}%
	\addtocounter{footnote}{-1}%
	\endgroup
}

\makeatletter
\patchcmd{\@maketitle}{\@title}{\fontsize{12}{14.4}\selectfont\@title}{}{}
\makeatother

\author{\textbf{{\fontfamily{ptm}\selectfont S. DELIKARAOGLOU, P. PINSON, R. ERIKSSON, T. WECKESSER}}\vspace{-0.35cm}}
\affil{\textbf{{\fontfamily{ptm}\selectfont Dept. of Electrical Engineering}} \\ \textbf{{\fontfamily{ptm}\selectfont Technical University of Denmark}} \\ \textbf{{\fontfamily{ptm}\selectfont Kgs. Lyngby 2800}} }

\begin{document}

\title{\textbf{Optimal Dynamic Capacity Allocation of HVDC Interconnections for Cross-border Exchange of Balancing Services in Presence of Uncertainty - Extended Version}\vspace{-0.3cm}}
\date{}
\maketitle
\thispagestyle{empty}

\vspace{1.7cm}

\section*{SUMMARY}

The deployment of large shares of stochastic renewable energy, e.g., wind power, may bring important economic and environmental benefits to the power system. Nonetheless, their efficient integration depends on the ability of the power system to cope with their inherent variability and the uncertainty arising from their partial predictability. Considering that the existing setup of the European electricity markets promotes the spatial coordination of neighbouring power systems only on the day-ahead market stage, regional system operators have to rely mainly on their internal balancing resources in order to guarantee system security. However, as power systems are forced to operate closer to their technical limits, where flexible generation resources become scarce, the conventional market paradigm may not be able to respond effectively on the wide range of uncertainty. 

The operational flexibility of the power system depends both on the technical parameters of its components, i.e., generators and transmission infrastructure, as well as on the operational practices that make optimal use of the available assets. This work focuses on alternative market designs that enable sharing of cross-border balancing resources between adjacent power systems through High Voltage Direct Current (HVDC) interconnections which provide increased controllability. In this context, we formulate a stochastic market-clearing algorithm that attains full spatio-temporal integration of reserve capacity, day-ahead and balancing markets. Against this benchmark we compare two deterministic market designs with varying degrees of coordination between the reserve capacity and energy services, both followed by a real-time mechanism. Our study reveals the inefficiency of deterministic approaches as the shares of wind power increase. Nevertheless, enforcing a tighter coordination between the reserves and energy trading floors may improve considerably the expected system cost compared to a sequential market design. Aiming to provide some insights for improvement of the sequential market-clearing, we analyse the effect of explicit transmission allocation between energy and reserves for different HVDC capacities and identify the market dynamics that dictate the optimal ratio.

\blfootnote{E-mail: stde@elektro.dtu.dk}

\section*{KEYWORDS}
Electricity markets, high voltage direct current (HVDC), reserves, transmission capacity allocation, stochastic programming, wind power.

\section{INTRODUCTION}
\vspace{-0.2cm}
The European power system is undergoing significant restructuring both in terms of generation capabilities and market architecture. The high shares of variable and partially predictable generation, e.g., wind power, pose new challenges to system operators due to the increased reserve requirements arising from the associated forecast errors. Now, although the appeal for improved efficiency and liquidity has already motivated the establishment of a single European day-ahead electricity market, the provision and the activation of reserves remain mainly on a national level. However, in this new operation paradigm, the integration of the reserve capacity and balancing markets seems to be a necessary step in order to enable the cross-border exchange of regulating services and reduce the balancing costs.

A major prerequisite for the deployment of cross-border reserves during real-time operation is the availability of transmission capacity between different areas of the power system. Meanwhile, the transmission grid should have adequate operational flexibility in order to allow for frequent re-dispatching of control reserves depending on the actual imbalances. HVDC interconnections can significantly contribute to a reliable solution owing to their increased controllability features that enable dynamic changes of the power flow direction which can be used to mitigate AC cross-border congestions \cite{chatzivasileiadis2013power}.

Given the current European market setup, where reserve capacity, forward and balancing market-clearings are performed sequentially, the ex-ante allocation of transmission capacity between these trading floors is necessary. However, this procedure is subject to uncertainty due to the significant time interval between day-ahead scheduling and actual operation. The two main approaches to cope with this uncertainty are stochastic programming and deterministic models \cite{morales2013integrating}. Despite the theoretical advances of the former, current operational practices are based on deterministic and usually static reserve requirements. For instance, in the new HVDC interconnection between Denmark and Norway (Skagerrak 4), 15\% of its capacity will be permanently reserved for exchange of ancillary services \cite{Energinet}. 


Previous studies \cite{gebrekirosoptimal}, \cite{farahmand2012balancing} have shown the benefit of cross-border reserve procurement in the Northern European power system based on deterministic reserve requirements. In this work, we model and assess the efficiency of transmission capacity allocation for the exchange of balancing power in systems with a large stochastic wind power component. In this context, we define two deterministic market setups which are compared with the ideal stochastic clearing. In addition, we analyse the sensitivity of the deterministic solution with respect to particular parameters, i.e., wind power penetration and shares of transmission capacity between different markets. The proposed methodology is applied on a modified power system based on the IEEE RTS-96 with additional HVDC interconnections.

\vspace{-0.3cm}


\section{MATHEMATICAL FORMULATION}
\label{sec:MathFor}
\vspace{-0.2cm}
This section provides the mathematical formulation of three alternative market-clearing models employed to settle the \textit{reserve capacity}, \textit{day-ahead} and \textit{balancing} markets. All the market setups considered in the present work are based on the following assumptions: \begin{inparaenum}[\itshape i\upshape)]
	\item only wind power uncertainty is considered assuming that load and availability of conventional generators are perfectly known,  
	\item demand is completely inelastic and
	\item the market settles on an hourly basis and each trading period is independent since no inter-temporal constraints, e.g., ramping limits, are considered.
\end{inparaenum} The market-clearing models presented in the remaining of this section are described in a compact form, while the full model formulation is available in the Appendix.  
\vspace{-0.35cm}
\subsection{Stochastic energy and reserve co-optimization} 
In order to determine the optimal amount of reserves and the day-ahead dispatch taking into account the anticipated real-time operation of the power system, market operator solves model \ref{StochMC} having a finite set of scenarios that captures the spatial dependence structure of the forecast errors. This formulation performs an intrinsic trade-off between the value of the available resources in the day-ahead (energy and reserve capacity) and the balancing (reserve deployment) stages, while it enables an implicit allocation of the transmission capacity among energy and reserve services.

\begin{subequations}
	\label{StochMC}
	\begin{align}
	\centering
	\label{eqn:S_ObjFun}
	&\underset{\Xi_{S}}{\text{Minimize}} \quad   C^D(p_C,p_W,\hat{r}) + \mathbb{E}_{\omega}\left[{C^B(\tilde{r}_{\omega})}\right]\\ 	
	&\text{subject to} \nonumber
	\end{align}
\begin{minipage}[t]{0.4\textwidth}
	\vspace{-24pt}
	\begin{align}
	 & h^D(p_C,p_W,\hat{\delta},\hat{f}) - l = 0 \label{eqn:S_Eq_DA} \\
	 & g_1^D(p_C,\hat{r}) \leq 0 \label{eqn:S_Ineq_DA}\\
	 & p_W \leq \overline{p}_W   \label{eqn:S_WindSchedule} \\
	 & g_2^D(\hat{\delta}, \hat{f}) \leq \overline{f}   \label{eqn:S_Ineq_DA_delta}\\
	 & g_3^D(\hat{r}) \leq \overline{R}   \label{eqn:S_Ineq_DA_Rmax}
	\end{align}		
\end{minipage}
\begin{minipage}[t]{0.6\textwidth}
		\vspace{-24pt}	
		\begin{align}
		 & h^B(\tilde{r}_\omega,\hat{\delta},\tilde{\delta}_\omega, \hat{f}, \tilde{f}_{\omega}) + \mathcal{W}_{\omega} - p_W = 0 , \forall \omega \in \Omega \label{eqn:S_Eq_Bal} \\
		 & g_1^B(\tilde{r}_\omega,\hat{r}) \leq 0, \; \forall \omega \in \Omega \label{eqn:S_Ineq_Bal2} \\
		 & g_2^B(\tilde{\delta}_\omega, \tilde{f}_{\omega}) \leq \overline{f}, \;\forall \omega \in \Omega \label{eqn:S_Ineq_Bal_delta}
		\end{align}	
\end{minipage}
\end{subequations}
\vspace{12pt}

where $\Xi_{S} = \{p_C, p_W,\hat{r}, \hat{\delta}, \hat{f}, \tilde{r}_\omega,\tilde{\delta}_\omega, \tilde{f}_{\omega}\}$ is the set of optimization variables and $\mathbb{E}_{\omega}[\cdot]$ the expectation operator over the scenario set $\Omega$ which describes wind power uncertainty. This is a two-stage stochastic programming problem, where vectors $p_C$ and $p_W$ denote the day-ahead dispatch of conventional and wind power generators, respectively. Vector $\hat{r}$ contains the reserve capacity scheduled for units that provide balancing services. All re-dispatching actions, i.e., up- and down-regulation, wind spillage and load shedding, that take place in real-time operation are included in vector $\tilde{r}_{\omega}$. The nodal voltage angles corresponding to day-ahead and balancing stage are denoted as $\hat{\delta}$ and $\tilde{\delta}_{\omega}$, while $\mathcal{W}_{\omega}$ corresponds to the realization of the random variable $\mathcal{W}$, modelling the actual wind power production in scenario $\omega$. The power flows on the HVDC links, denoted as $\hat{f}$ at the day-ahead stage and $\tilde{f}_{\omega}$ during the balancing stage, are considered fully controllable and any voltage limitations and active power losses are neglected. Constants $l$, $\overline{f}$ and $\overline{p}_W$ denote the load at each node of the system, the transmission capacity limits and the wind power installed capacity, respectively. 

Constraints \eqref{eqn:S_Eq_DA} and \eqref{eqn:S_Eq_Bal} are the nodal power balance equations in the day-ahead and balancing stages. The day-ahead energy dispatch and reserve provision are limited by the generator maximum production through constraints \eqref{eqn:S_Ineq_DA} and \eqref{eqn:S_WindSchedule}. Similarly, constraints \eqref{eqn:S_Ineq_DA_Rmax} ensure that the reserved capacity of each generator does not exceed its reserve quantity offer $\overline{R}$, while constraints \eqref{eqn:S_Ineq_Bal2} guarantee that the balancing actions are limited by the corresponding reserved capacity or rational constraints, i.e., load shedding and wind spillage are lower than or equal to the load and the realized wind power production, respectively. Finally, the set of constraints \eqref{eqn:S_Ineq_DA_delta} and \eqref{eqn:S_Ineq_Bal_delta} enforces the transmission capacity limits at the day-ahead and balancing stages. 

\vspace{-0.2cm}

\subsection{Deterministic market-clearing}


In order to avoid the computational burden of the stochastic market-clearing approach, market operators employ deterministic models which have to procure reserve capacity based on explicit reserve requirements in order to deal with unforeseen events during real-time operation. This section provides the mathematical formulation of two deterministic approaches; the first allows for simultaneous procurement of energy and reserves while the second decouples completely these two services through a sequential settlement of the energy and reserve capacity markets. Both market organizations include a real-time balancing mechanism to compensate any deviations from the day-ahead schedule.

\vspace{-0.2cm}
 
\subsubsection{Energy and reserve capacity co-optimization}

The day-ahead market-clearing problem for energy and reserve co-optimization is described through model \ref{ConvMC_DA}, followed by the balancing market settlement of model \ref{ConvMC_B}. This market setup finds the optimal day-ahead dispatch (denoted by the superscript "$*$") and yields an implicit optimal transmission capacity allocation between energy and reserves, based on a single-valued forecast $\widehat{\mathcal{W}}$ (conditional expectation) of wind power production, taking into account the explicit system-wide reserve requirements $\overline{RR}_s$. The optimization variables of the day-ahead and the balancing market-clearing problems are included in the sets $\Xi_{C}^{D} = \{p_C, p_W, \hat{r}, \hat{\delta}, \hat{f} \}$ and $\Xi_{C}^{B} = \{\tilde{r}, \tilde{\delta}, \tilde{f} \}$, respectively. 

\begin{minipage}[t]{0.4\textwidth}
	\begin{center}
	Day-ahead market-clearing
	\end{center}
	\vspace{-16pt}
\begin{subequations}
	\label{ConvMC_DA}
	\begin{align}
	\label{eqn:DA_ObjFun}
	&\underset{\Xi_{C}^{D}}{\text{Minimize}} \quad C^D(p_C,p_W,\hat{r})  \\ 
	&\text{subject to} \nonumber \\
	& h^D(p_C,p_W,\hat{\delta},\hat{f}) - l = 0  \label{eqn:C_Eq_DA} \\
	& g_1^D(p_C,\hat{r}) \leq 0 \label{eqn:C_Ineq_PR}\\
	& p_W \leq \widehat{\mathcal{W}}   \label{eqn:C_WindSchedule} \\
	& g_2^D(\hat{\delta}, \hat{f}) \leq \overline{f}   \label{eqn:C_Ineq_DA_delta}\\	
	& g_3^D(\hat{r}) \leq \overline{R}   \label{eqn:C_Ineq_DA_Rmax}\\	
	& g_4^D(\hat{r}) \geq \overline{RR}_s   \label{eqn:C_Ineq_DA_RR}
	\end{align}
\end{subequations}
\end{minipage}
\begin{minipage}[t]{0.54\textwidth}	
	\begin{center}
		Balancing market-clearing	
			\vspace{-16pt}
	\end{center}
\begin{subequations}
	\label{ConvMC_B}
	\begin{align}
	\label{eqn:B_ObjFun}
	&\underset{\Xi_{C}^{B}}{\text{Minimize}} \quad C^B(\tilde{r})\\
	&\text{subject to} \nonumber \\
	& h^B(\tilde{r},\hat{\delta}^{*},\tilde{\delta}, \hat{f}^{*}, \tilde{f})+  \mathcal{W} - p_W^{*} = 0 \label{eqn:C_Eq_Bal} \\
	& g_1^B(\tilde{r}) \leq \hat{r}^{*}  \label{eqn:C_Ineq_Bal1} \\
	& g_2^B(\tilde{\delta}, \tilde{f}) \leq \overline{f}, \label{eqn:C_Ineq_Bal2}
	\end{align}
\end{subequations}
\end{minipage}
\vspace{-0.2cm}
\subsubsection{Sequential clearing of reserve capacity and energy markets}
The separation of reserve capacity and energy markets is implemented using models \ref{ConvMC_ResCap} and \ref{ConvMC_DA_Area}. This market architecture requires the explicit allocation of transmission capacity for reserve exchange, denoted as $X$, as well as the definition of reserve requirements $\overline{RR}_{a}$ for each area $a$ of the power system. Model \ref{ConvMC_ResCap} determines the least-cost reserve capacity $\hat{r}_a^{*}$ which is used as input in the energy-only day-ahead clearing \ref{ConvMC_DA_Area} whose optimization variables are contained in the set $\Xi_{C}^{SD} = \{p_C, p_W, \hat{\delta}, \hat{f}\}$. Any corrective actions during real-time operation are compensated through a balancing market as of  model \ref{ConvMC_B} considering that a single pool of control resources is accessible from all areas.

\vspace{0.2cm}

\begin{minipage}[t]{0.36\textwidth}
	\begin{center}
	Reserve capacity market-clearing
		\vspace{-16pt}
	\end{center}
\begin{subequations}
	\label{ConvMC_ResCap}
	\begin{align}
	\label{eqn:RCap_ObjFun}
	&\underset{\hat{r}_a}{\text{Minimize}} \quad C^R(\hat{r}_{a})  \\ 
	&\text{subject to} \nonumber \\
	& g_1^R(\hat{r}_{a}) \leq \overline{R}   \label{eqn:R_Rmax}\\	
	& g_2^R(\hat{r}_{a}) \geq \overline{RR}_{a}   \label{eqn:R_RR_Area}\\	
	& g_3^R(\hat{r}_a) \leq X\overline{f}   \label{eqn:R_TransCap}
	\end{align}
\end{subequations}
\end{minipage}
\begin{minipage}[t]{0.54\textwidth}
	\begin{center}
		Day-ahead market-clearing
			\vspace{-16pt}
	\end{center}
	\begin{subequations}
		\label{ConvMC_DA_Area}
		\begin{align}
		\label{eqn:DA_ObjFun}
		&\underset{\Xi_{C}^{SD}}{\text{Minimize}} \quad  C^D(p_C,p_W)  \\ 
		&\text{subject to} \nonumber \\ 
		& h^D(p_C,p_W,\hat{\delta},\hat{f}) - l = 0  \label{eqn:C_Eq_DA} \\
		& g_1^D(p_C,\hat{r}^{*}_{a}) \leq 0 \label{eqn:C_Ineq_PR}\\
		& p_W \leq \widehat{\mathcal{W}}   \label{eqn:C_WindSchedule} \\
		& g_2^D(\hat{\delta}, \hat{f}) \leq (1-X)\overline{f}   \label{eqn:C_Ineq_DA_delta}		
		\end{align}
	\end{subequations}
\end{minipage}

\vspace{-0.3cm}

\section{RESERVE REQUIREMENTS AND WIND POWER UNCERTAINTY}
\label{sec:Uncrt}
\vspace{-0.25cm}
The reserve requirements are calculated using a probabilistic approach similar to \cite{Doherty2005}, aiming to cover a pre-specified interval $\xi$ of the wind power predictive distribution. Write $\bar{\alpha}, \underline{\alpha} \in [0,1]$ the appropriate nominal proportions of the predictive cumulative distribution function (CDF) $\hat{F}$ such that $\xi=\overline{\alpha}-\underline{\alpha}$, e.g., for reliability level $\xi=99\%$ it follows that $\overline{\alpha}=0.995$ and $\underline{\alpha}=0.005$. Then, the reserve requirements for up- and down-regulation from wind farm $s$, denoted as $RR^{+}_s$ and $RR^{-}_s$ respectively, are calculated as
\begin{align}
RR^{+}_s = \left(\hat{\mu}_s - \hat{F}^{-1}_s (\underline{\alpha})\right) \bar{p}_{W}^{s}, \quad
RR^{-}_s =\left(\hat{F}^{-1}_s (\overline{\alpha}) - \hat{\mu}_s\right) \bar{p}_{W}^{s}
\label{eqn:ResReq}
\end{align}
where $\hat{\mu}_{s}$ is the mean of $\hat{F}_{s}$. System-wide reserve requirements are obtained from \eqref{eqn:ResReq}, replacing $\hat{F}_s$ with the CDF of the whole wind power portfolio and $\bar{p}_{W}^{s}$ with the total installed wind power capacity. 


The stochastic output of wind farm $s$, normalized by its nominal capacity, is modelled using a Beta distribution $B_s(\alpha_s, \beta_s)$ as proposed in \cite{Fabbri2005}. To obtain a set $\Omega$ of spatially correlated scenarios we follow the methodology presented in \cite{Sdeli}.  Denote as $X^s_{\omega}$ the realization drawn from a multivariate Normal distribution $\mathcal{N}(0, \Sigma)$, where $\Sigma$ the covariance matrix with diagonal elements equal to 1 and off-diagonal elements equal to the spatial correlation coefficient $\rho$ of the corresponding wind locations. Then, the $\mathcal{W}^s_{\omega}$ wind power production scenario is generated using the following transformation:
\begin{equation}
Y^s_{\omega} = \Phi (X^s_{\omega}) \longrightarrow \mathcal{W}^s_{\omega} = \hat{F}_{B|s}^{-1}\left(Y^s_{\omega}\right)
\end{equation}
where $\Phi$ and $\hat{F}_{B|s}$ denote the Gaussian and Beta CDF, respectively.

\section{CASE STUDY}
\vspace{-0.2cm}
In this section we evaluate the performance of the market-clearing algorithms described in section \ref{sec:MathFor}, using a power system comprised of areas 1 and 2 of the IEEE RTS-96 system presented in \cite{REALTestSystems}. The offer price of each conventional generator in the day-ahead and balancing markets is set equal to the average of the three marginal cost segments provided in \cite{REALTestSystems}, while wind power offers its energy production at zero price. In addition, the minimum power production of all generators is assumed to be zero and the value of lost load is set to 1000\$/MWh. OCGT plants provide their full capacity in the reserve market at cost per MW  equal to 25\% of their day-ahead offer price. The IGCC and CCGT units offer 40\% and 25\% of their capacity for reserve provision at a price equal to 10\% and 5\% of their marginal cost. Considering their limited flexibility, nuclear and coal units do not provide any balancing services. In order to focus our analysis on the HVDC tie-line, we follow a zonal pricing scheme similar to the NordPool market, neglecting the intra-area network constraints. This allows to model wind power production in areas 1 and 2 by employing two Beta distributions with shape parameters $(\alpha, \beta)$ equal to (3.78, 1.62) and (5.67, 6.48), respectively. A set of 100 equiprobable scenarios is obtained assuming a correlation coefficient $\rho$ equal to 0.35.

\begin{figure}
	\centering
	\begin{minipage}{.5\textwidth}
		\centering
		\includegraphics[width=1\textwidth]{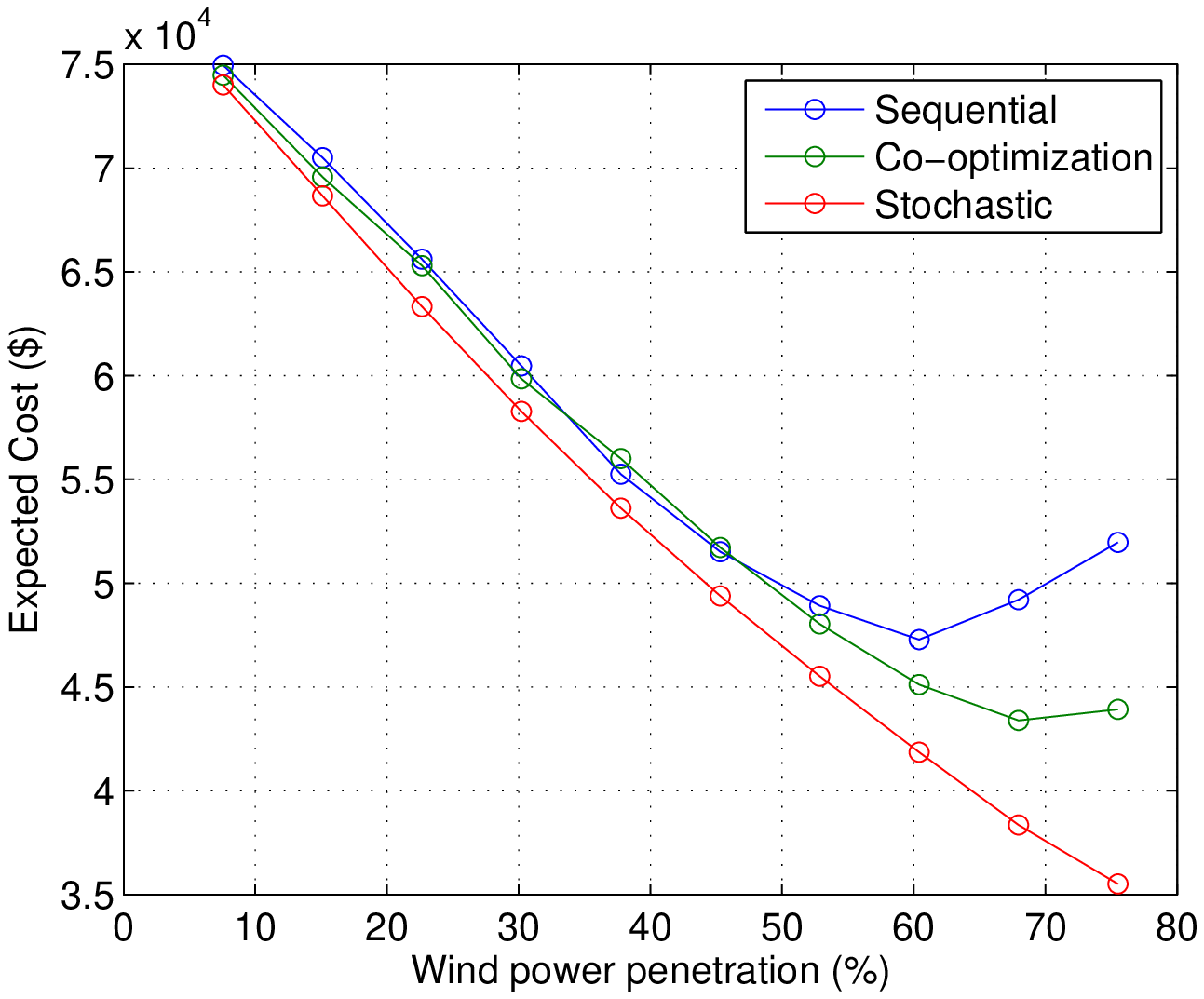}
		\vspace{-0.7cm}
		\captionof{figure}{Expected cost as a function of wind power penetration \\ for different market-clearing models ($X=15\%$).}
		\label{fig:CostComp}
	\end{minipage}%
	\begin{minipage}{.5\textwidth}
		\centering
		\includegraphics[width=1\textwidth]{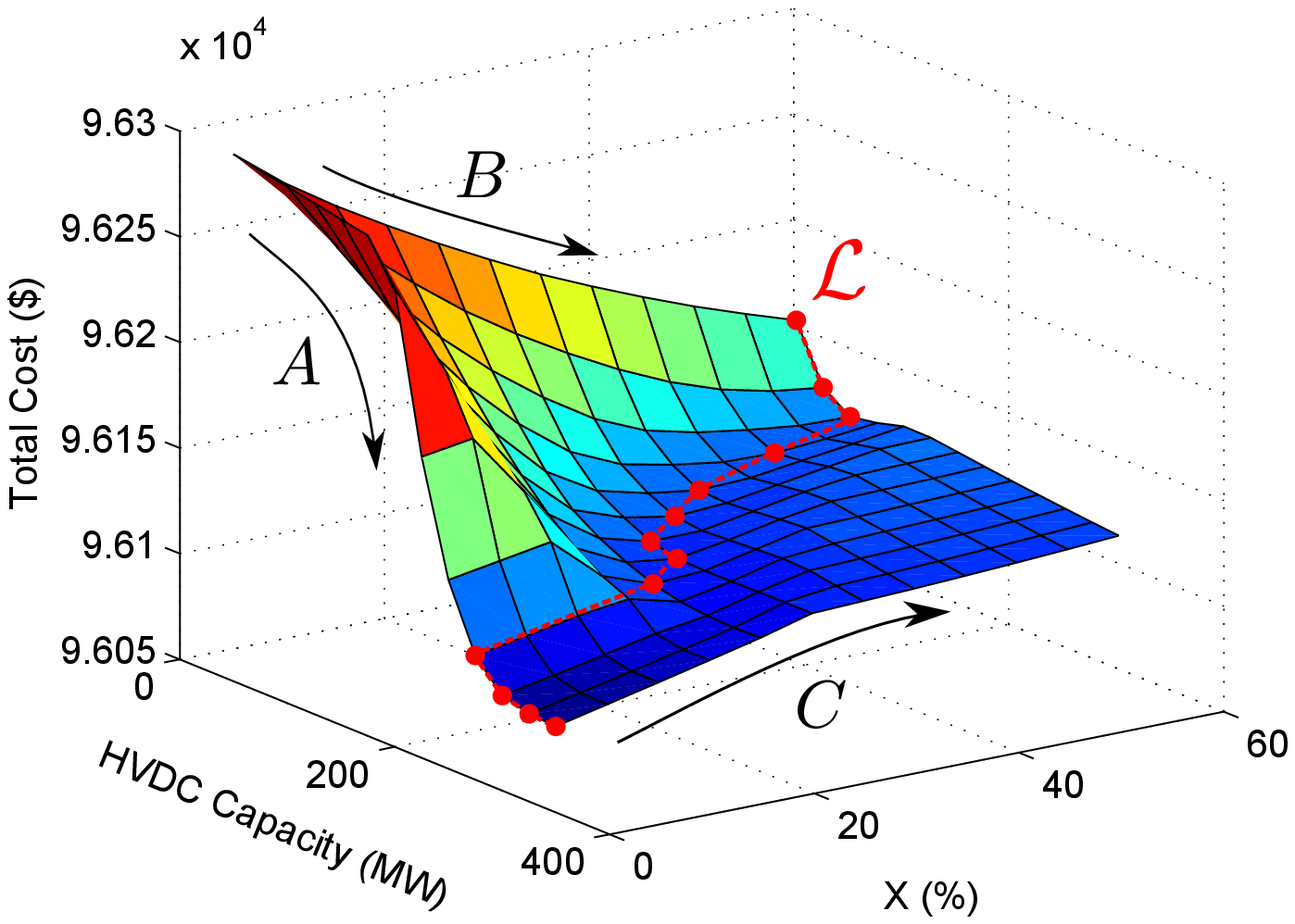}
		\vspace{-0.7cm}
		\captionof{figure}{Impact of HVDC capacity and transmission allocation $X$ on the expected cost of sequential market-clearing.}
		\label{fig:Cost3D}
	\end{minipage}
	\vspace{-0.4cm}
\end{figure}

Figure \ref{fig:CostComp} shows the evolution of the expected power system operation cost as wind penetration increases, for the three considered market setups and an HVDC tie-line capacity of 200 MW. The ratio of the installed wind power capacity between areas 1 and 2 is fixed to 2:1 and wind power penetration is defined as a percentage of total system load. It can be observed that stochastic market-clearing outperforms both deterministic approaches in the whole range of wind penetration levels, exploiting the advanced information on the spatial characteristics of wind forecast uncertainty and managing to fully capture the benefits of cost-free wind power. The performance of the two deterministic market setups is similar up to a penetration level of 45\%, while beyond this wind power share the deterministic energy and reserve co-optimization model achieves significantly lower expected costs. However, the inefficiency of both deterministic models becomes apparent for large penetration levels, where increasing installed wind power capacity leads to higher expected cost. It should be noted that the curves in figure \ref{fig:CostComp} are calculated using in-sample analysis, i.e., the same scenario set $\Omega$ describing wind uncertainty in the day-ahead stage is also used to calculate the expected balancing cost of each market design.

In order to study the effect of the explicit transmission allocation $X$ on the expected system cost, we consider a more penalizing reserve market where OCGT, IGCC and CCGT units provide 50\%, 40\% and 40\% of their capacity for reserve provision at a price equal to 15\%, 30\% and 30\% of their day-ahead offer, respectively. For a constant wind penetration level of 24\%, we perform a 'grid-search' on the HVDC tie-line capacity and the parameter $X$, to obtain the corresponding expected system cost displayed in figure \ref{fig:Cost3D}. The shape of this surface allows to identify three main directions aiming to analyse the effect of these parameters on the overall market efficiency. Note that the dotted line $\mathcal{L}$ represents the locus of the least-cost points for different HVDC capacities. Moving along direction $A$ on the left-hand side of $\mathcal{L}$, the expected cost reduces as the tie-line capacity increases since a larger pool of common resources is accessible from both areas. A similar trend applies for larger values of $X$ (direction $B$), where the reserve capacity market does not alter significantly the optimal day-ahead settlement by imposing additional constraints. For example, in a segmented market of reserve resources (low $X$), expensive units may be assigned to provide  downward reserve, enforcing a corresponding lower bound on the energy dispatch, out of the merit order. Large values of $X$ on high HVDC capacity (direction $C$), allow the reserve market to pick more economical resources based solely on their capacity payment. However, this translates into increased expected balancing cost if mainly excess production situations occur during actual operation, i.e., generators with low marginal (and reserve capacity) cost are willing to buy back their production in lower price.

\vspace{-0.5cm}

\section{CONCLUSIONS AND FUTURE WORK}
\vspace{-0.2cm}
This paper considers the problem of optimal allocation of flexible HVDC transmission for the cross-border exchange of balancing reserves in presence of uncertain wind power generation. Setting as benchmark the stochastic energy and reserves co-optimization, we analyse the efficiency of two deterministic market designs to accommodate large shares of wind power and we show that both versions behave fairly similar, only truly revealing the benefits from co-optimization when getting close to the system limits. Despite that the relative difference in market efficiencies may be highly sensitive to the system properties, the ranking is expected to stay the same. A wide range of plausible system setups is explored in order to perceive the market dynamics defining the optimal transmission allocation in the sequential market-clearing. It can be concluded that there generally exists an optimal allocation to be made, which however may dynamically vary depending on generation, load and system uncertainties.  

Future work may focus on the reformulation of the deterministic energy and reserve co-optimization algorithm by defining per area reserve requirements in order to reveal the optimal transmission capacity allocation of this model. This could serve as an insight towards an explicit rule that determines the optimal percentage of the tie-line capacity to be designated to reserve exchange in the sequential design. Further extensions could investigate the issues related to HVDC and AC network coordination as well as integrate the aspect of cross-border reserve provision in the HVDC expansion planning.

  

\vspace{-0.4cm}

\section*{BIBLIOGRAPHY}
\vspace{-.3cm}

\bibliographystyle{IEEEtran}
\bibliography{IEEEabrv,CIGRE}

{\lmr

\section*{APPENDIX}



\section*{Nomenclature}
\subsection*{Indices and Sets}


\begin{table}[h]
	\renewcommand{\familydefault}{}\normalfont
	\begin{tabular}{ll}
		$n$               & Index of system buses.                      \\
		$i$               & Index of dispatchable units.                \\
		$j$               & Index of wind power units.                  \\
		$l$               & Index of transmission lines.                \\
		$\omega$          & Index of wind power scenarios.              \\
		$L^{\text{AC}}$   & Set of AC transmission lines.               \\
		$L^{\text{DC}}$   & Set of HVDC transmission lines.             \\
		$\mathcal{M}_n^J$ & Set of wind power units located at bus $n$. \\
		$\mathcal{M}_n^I$ & Set of dispatchable units located at bus $n$.                                           
	\end{tabular}
\end{table}

\subsection*{Parameters}


\begin{table}[h]
	\renewcommand{\familydefault}{}\normalfont
	\begin{tabular}{ll}
		$A_{ln}$                                & \begin{tabular}[c]{@{}l@{}}Line-to-bus incidence matrix. The elements of the incidence matrix $A_{ln}$ \\ are equal to -1/1 if line $l$ leaves/enters bus $n$ and 0 otherwise.\end{tabular} \\
		$D_n$                                   & Demand at node $n$ {[}MW{]}.                                                                                                                     \\
		$C_{i}$                                 & Day-ahead offer price of unit $i$ {[}\$/MWh{]}.                                                                                                  \\
		$C^{\text{R}+}_i,C^{\text{R}-}_i$       & Up/down reserve capacity offer price of unit $i$ {[}\$/MW{]}.                                                                                    \\
		$C^{\text{sh}}$                         & Value of involuntarily shed load {[}\$/MWh{]}.                                                                                                   \\
		$P^{\text{max}}_i,P^{\text{W,max}}_j$   & Capacity of units $i$ and $j$ {[}MW{]}.                                                                                                          \\
		$\bar{P}^{\text{W}}_{j}$                & Forecast mean of wind power production {[}MW{]}.                                                                                                 \\
		$P^*_{j\omega}$                         & Wind power realization in scenario $\omega$ {[}MW{]}.                                                                                            \\
		$R^{+,\text{max}}_i,R^{-,\text{max}}_i$ & Maximum up/down reserve provided by unit $i$ {[}MW{]}.                                                                                           \\
		$RR^{\text{s},+},RR^{\text{s},-}$       & Maximum up/down system reserve requirements {[}MW{]}.                                                                                            \\
		$f^{\text{max}}_{l}$                    & Maximum capacity of line $l$ {[}MW{]}.                                                                                                           \\
		$B_l$                                   & Absolute value of the susceptance of AC line $l$ {[}per unit{]}.                                                                                
	\end{tabular}
\end{table}

\subsection*{Variables}


\begin{table}[h]
	\renewcommand{\familydefault}{}\normalfont
	\begin{tabular}{ll}
		$p_i,p_{j}^\text{W}$                                               & Day-ahead dispatch of units $i$ and $j$ {[}MW{]}.                      \\
		$R^{+}_i,R^{-}_i$                                                  & Up/down reserve capacity scheduled for unit $i$ {[}MW{]}.              \\
		$r^+_{i\omega},r^-_{i\omega}$                                      & Up/down regulation provided by unit $i$ in scenario $\omega$ {[}MW{]}. \\
		$p^{\text{spill}}_{j\omega}$                                       & Wind power spilled by unit $j$ in scenario $\omega$,{[}MW{]}.          \\
		$l^{\text{sh}}_{n\omega}$                                          & Load shedding at node $n$ in scenario $\omega$ {[}MW{]}.               \\
		$\hat{\delta}_n$                                                   & Voltage angle at node $n$ at the day-ahead market stage {[}rad{]}.     \\
		$\tilde{\delta}_{n\omega}$                                         & Voltage angle at node $n$ in scenario $\omega$ {[}rad{]}.              \\
		$\hat{f}^{\text{AC}}_l, \hat{f}^{\text{DC}}_l$                     & Power flow in AC/HVDC line $l$ at day-ahead market stage {[}MW{]}.     \\
		$\tilde{f}^{\text{AC}}_{l\omega}, \tilde{f}^{\text{DC}}_{l\omega}$ & Power flow in AC/HVDC line $l$ in scenario $\omega$ {[}MW{]}.         
	\end{tabular}
\end{table}

\newpage
\section*{Stochastic energy and reserve co-optimization}






\begin{subequations}
\begin{flalign}
&\underset{\Xi^{S}}{\text{Min.}} \thinspace \sum_{i\in I}\left({C}_{i}{p}_{i} + C^{\text{R}+}_i R^{+}_{i} + C^{\text{R}-}_i R^{-}_{i}\right) + \sum_{\omega \in \Omega} \pi_{\omega} \left[\sum_{i \in I}{C}_{i}\left({r}_{i\omega}^{+}-{r}_{i\omega}^{-}\right)+\sum_{n \in N} {C}^{\text{sh}}l^{\text{sh}}_{n\omega}\right]\label{eq:low-obj} \\
& \text{subject to} \nonumber \\
& \sum_{j \in \mathcal{M}_n^J} p_j^\text{W} + \sum_{i\in \mathcal{M}_n^I}{p}_{i} - D_n - \sum_{l \in L^{AC}\cup L^{DC}}  A_{ln} \left( \hat{f}^{\text{AC}}_{l} + \hat{f}_l^{DC} \right) = 0, \quad \forall n  \label{eq:SP-da-bal-constr} \\
& \sum_{i \in \mathcal{M}_n^I} \left(r_{i\omega}^{+}-r_{i\omega}^{-} \right) + l^{\text{sh}}_{n\omega} + \sum_{j \in \mathcal{M}_n^J}\left(P_{j\omega}^*-p_j^\text{W}-p^{\text{spill}}_{j\omega}\right) \nonumber \\
& + \sum_{l \in L^{AC}\cup L^{DC}}  A_{ln} \left(\hat{f}^{\text{AC}}_{l} -\tilde{f}^{\text{AC}}_{l\omega} + \hat{f}^{DC}_{l} - \tilde{f}^{DC}_{l\omega} \right) = 0, \quad  \forall n, \forall \omega  \label{eq:SP-rt-bal-constr} \\
& p_j^\text{W} \leq P^{\text{W,max}}_{j},  \quad \forall j  \label{eq:low-da-wind-max}\\
& p_i + R^+_{i} \leq P^{\text{max}}_{i}, \quad \forall i  \label{eq:SP-rt--conv-max}\\
& p_i - R^-_{i} \geq 0, \quad \forall i  \label{eq:SP-rt--conv-min}\\
& R^+_{i} \le R^{+,\text{max}}_{i}, \quad \forall i \label{eq:SP-rescap-up-max}\\
& R^-_{i} \le R^{-,\text{max}}_{i}, \quad \forall i \label{eq:SP-rescap-dn-max}\\
& r^+_{i\omega} \leq R^+_i, \quad \forall i, \forall \omega  \label{eq:SP-rt-conv-res-up-max}\\
& r^-_{i\omega} \leq R^-_i, \quad \forall i, \forall \omega  \label{eq:SP-rt-conv-res-dn-max}\\
& \hat{f}^{\text{AC}}_{l} = B_l \sum_n A_{ln} \hat{\delta}_n, \quad \forall l \in L^{AC} \label{eq:SP-da-flowAC-constr}\\
& \tilde{f}^{\text{AC}}_{l\omega} = B_l \sum_n A_{ln} \tilde{\delta}_{n\omega}, \quad \forall l \in L^{AC}, \forall \omega \label{eq:SP-rt-flowAC-constr}\\
& -f^{\text{max}}_{l} \le \hat{f}^{\text{AC}}_{l} \le f^{\text{max}}_{l}, \quad \forall l \in L^{AC} \label{eq:SP-da-flowAC-max}\\
& - f^{\text{max}}_{l} \le \tilde{f}^{\text{AC}}_{l\omega} \le f^{\text{max}}_{l}, \quad \forall l \in L^{AC}, \forall \omega \label{eq:SP-rt-flowAC-max}\\
& -f^{\text{max}}_{l} \le \hat{f}^{DC}_l \le f^{\text{max}}_{l}, \quad \forall l \in L^{DC} \label{eq:SP-da-flowDC-constr}\\
& -f^{\text{max}}_{l} \le \tilde{f}^{DC}_{l\omega} \le f^{\text{max}}_{l}, \quad \forall l \in L^{DC}, \forall \omega \label{eq:SP-rt-flowDC-constr}\\
& l^{\text{sh}}_{n\omega} \leq D_n, \quad \forall n, \forall \omega  \label{eq:SP-rt-load-shed-max}\\
& p^{\text{spill}}_{j\omega} \leq P^*_{j\omega}, \quad \forall j, \forall \omega  \label{eq:SP-rt-wind-spill-max}\\
& \hat{\delta}_1 = 0; \enspace \tilde{\delta}_{1\omega} = 0, \quad \forall \omega \label{eq:SP-ref-node}\\
& p_i \geq 0 , \forall i; \enspace p_j^{\text{W}} \geq 0 , \forall j; \enspace R^{+}_{i},R^{-}_{i}\geq 0, \forall i;  \nonumber \\     
& r^+_{i\omega},r^-_{i\omega}\geq 0, \forall i, \forall \omega; \enspace p^{\text{spill}}_{j\omega} \geq 0, \forall j, \forall \omega; \enspace l^{\text{sh}}_{n\omega} \geq 0, \forall n, \forall \omega;  \nonumber \\ 
& \hat{\delta}_n \thinspace \text{free}, \forall n; \enspace \tilde{\delta}_{n\omega} \thinspace \text{free}, \forall n, \forall \omega; \nonumber \\ 
& \hat{f}_l^{AC} \thinspace \text{free}, \forall l \in  L^{AC}; \enspace \tilde{f}_{l\omega}^{AC} \thinspace \text{free}, \forall l \in  L^{AC}, \forall \omega; \nonumber \\ 
& \hat{f}_l^{DC} \thinspace \text{free}, \forall l \in  L^{DC}; \enspace \tilde{f}_{l\omega}^{DC} \thinspace \text{free}, \forall l \in  L^{DC}, \forall \omega
\end{flalign}
\end{subequations}
The set of optimization variables \\
 $\Xi^{S}=\big\{p_i,p_{j}^\text{W},{R}^{+}_{i},{R}^{-}_{i},{r}^{+}_{i\omega},{r}^{-}_{i\omega}, p^{\text{spill}}_{j\omega}, {l}^{\text{sh}}_{n\omega}, \hat{\delta}_n, \tilde{\delta}_{n\omega}, \hat{f}^{AC}_l, \tilde{f}^{AC}_{l\omega}, \hat{f}^{DC}_l, \tilde{f}^{DC}_{l\omega}, \enspace \forall i, j, n, l, \omega \big\}$. 

\section*{Deterministic energy and reserve capacity co-optimization}

\subsection*{Day-ahead market-clearing: Energy and reserve capacity}

\begin{subequations}
\begin{flalign}
&\underset{\Xi^{C}_{DA}}{\text{Min.}} \thinspace \sum_{i\in I}\left({C}_{i}{p}_{i} + C^{\text{R}+}_i R^{+}_{i} + C^{\text{R}-}_i R^{-}_{i}\right) \\
& \text{subject to} \nonumber \\	
& \sum_{j \in \mathcal{M}_n^J} p_j^\text{W} + \sum_{i\in \mathcal{M}_n^I}{p}_{i} - D_n - \sum_{l \in L^{AC}\cup L^{DC}}  A_{ln} \left( \hat{f}^{\text{AC}}_{l} + \hat{f}_l^{DC} \right) = 0, \quad \forall n  \label{eq:CV-da-bal-constr} \\
& p_j^W \leq \bar{P}^{\text{W}}_{j},  \quad \forall j  \label{eq:CV-da-wind-max}\\
& p_i + R^+_{i} \leq P^{\text{max}}_{i}, \quad \forall i  \label{eq:CV-rt--conv-max}\\
& p_i - R^-_{i} \geq 0, \quad \forall i  \label{eq:CV-rt--conv-min}\\
& R^+_{i} \le R^{+,\text{max}}_{i}, \quad \forall i \label{eq:CV-rescap-up-max}\\
& R^-_{i} \le R^{-,\text{max}}_{i}, \quad \forall i \label{eq:CV-rescap-dn-max}\\
& \sum_{i} R_i^{+} \geq RR^{\text{s},+} \label{eq:Res-req-up}\\
& \sum_{i} R_i^{-} \geq RR^{\text{s},-} \label{eq:Res-req-dn}\\
& \hat{f}^{\text{AC}}_{l} = B_l \sum_n A_{ln} \hat{\delta}_n, \quad \forall l \in L^{AC} \label{eq:CV-da-flowAC-constr}\\
& -f^{\text{max}}_{l} \le \hat{f}^{\text{AC}}_{l} \le f^{\text{max}}_{l}, \quad \forall l \in L^{AC} \label{eq:CV-da-flowAC-max}\\			
& -f^{\text{max}}_{l} \le \hat{f}^{DC}_l \le f^{\text{max}}_{l}, \quad \forall l \in L^{DC} \label{eq:CV-da-flowDC-constr}\\
& \hat{\delta}_1 = 0\\
& p_i \geq 0 , \forall i; \enspace p_j^{\text{W}} \geq 0 , \forall j; \enspace R^{+}_{i},R^{-}_{i}\geq 0, \forall i; \nonumber \\ 
& \hat{\delta}_n \thinspace \text{free}, \forall n; \hat{f}_l^{AC} \thinspace \text{free}, \forall l \in  L^{AC}; \hat{f}_l^{DC} \thinspace \text{free}, \forall l \in  L^{DC}
\end{flalign}
\end{subequations}
The set of optimization variables $\Xi^{C}_{DA}=\big\{p_i,p_{j}^\text{W},{R}^{+}_{i},{R}^{-}_{i}, \hat{\delta}_n, \hat{f}^{AC}_l, \hat{f}^{DC}_l,  \enspace \forall i, j, n, l \big\}$.

\subsection*{Balancing market-clearing}

\begin{subequations}
	\begin{flalign}
	&\underset{\Xi^{C}_{B}}{\text{Min.}} \thinspace \sum_{i \in I}{C}_{i}\left({r}_{i}^{+}-{r}_{i}^{-}\right)+\sum_{n \in N} {C}^{\text{sh}}l^{\text{sh}}_{n} \\
	& \text{subject to} \nonumber \\
	& \sum_{i \in \mathcal{M}_n^I} \left(r_{i}^{+}-r_{i}^{-} \right) + l^{\text{sh}}_n + \sum_{j \in \mathcal{M}_n^J}\left(P_{j}^*-p^\text{W}_j-p^{\text{spill}}_{j}\right) \nonumber \\
	& + \sum_{l \in L^{AC}\cup L^{DC}}  A_{ln} \left(\hat{f}^{\text{AC}}_{l} -\tilde{f}^{\text{AC}}_{l} + \hat{f}^{DC}_{l} - \tilde{f}^{DC}_{l} \right) = 0, \quad  \forall n   \label{eq:CV-rt-bal-constr} \\	
	& r^+_{i} \leq R^+_i, \quad \forall i  \label{eq:CV-rt-conv-res-up-max}\\
	& r^-_{i} \leq R^-_i, \quad \forall i  \label{eq:CV-rt-conv-res-dn-max}\\
	& \tilde{f}^{\text{AC}}_{l} = B_l \sum_n A_{ln} \tilde{\delta}_{n}, \quad \forall l \in L^{AC}  \label{eq:CV-rt-flowAC-constr}\\
	& - f^{\text{max}}_{l} \le \tilde{f}^{\text{AC}}_{l} \le f^{\text{max}}_{l}, \quad \forall l \in L^{AC}  \label{eq:CV-rt-flowAC-max}\\
	& -f^{\text{max}}_{l} \le \tilde{f}^{DC}_{l} \le f^{\text{max}}_{l}, \quad \forall l \in L^{DC}  \label{eq:CV-rt-flowDC-constr}\\
	& l^{\text{sh}}_{n} \leq D_n, \quad \forall n   \label{eq:CV-rt-load-shed-max}\\
	& p^{\text{spill}}_{j} \leq P^*_{j}, \quad \forall j  \label{eq:CV-rt-wind-spill-max}\\	
	& \tilde{\delta}_{1} = 0 \label{eq:CV-ref-node}\\
	& r^+_{i},r^-_{i}\geq 0, \forall i; \enspace p^{\text{spill}}_{j} \geq 0, \forall j; \enspace l^{\text{sh}}_{n}, \forall n \geq 0; \nonumber \\ 
	& \tilde{\delta}_{n} \thinspace \text{free}, \forall n; \tilde{f}_{l}^{AC} \thinspace \text{free}, \forall l \in  L^{AC}; \enspace \tilde{f}_{l}^{DC} \thinspace \text{free}, \forall l \in  L^{DC}		
	\end{flalign}
\end{subequations}
The set of optimization variables $\Xi^{C}_{B}=\big\{{r}^{+}_{i},{r}^{-}_{i}, p^{\text{spill}}_{j} , {l}^{\text{sh}}_{n}, \tilde{\delta}_{n}, \tilde{f}^{AC}_{l}, \tilde{f}^{DC}_{l}, \enspace \forall i, j, n, l  \big\}$.

\newpage
\section*{Sequential clearing of reserve capacity and energy markets}

Additional nomenclature:

\begin{table}[h]
	\renewcommand{\familydefault}{}\normalfont
	\begin{tabular}{ll}
		$a, b$                   & Index of system areas.                                                                                                               \\
		$\mathcal{M}_a^I$        & Set of units located in area $a$.                                                                                                    \\
		$\mathcal{M}_{a b}^L$    & Set of lines connecting area $a$ with area $b$.                                                                                      \\
		$X_l$                    & Transmission capacity of line $l$ allocated to reserve exchange {[}\%{]}.                                                            \\
		$RR_a^{+}, RR_a^{-}$     & Up/down reserve requirements in area $a$ {[}MW{]}.                                                                                   \\
		$R^{+}_{ia}, R^{-}_{ia}$ & \begin{tabular}[c]{@{}l@{}}Up/down reserve capacity provided by unit $i$ \\ to cover requirements of area $a$ {[}MW{]}.\end{tabular}
	\end{tabular}
\end{table}


\subsection*{Reserve capacity market-clearing}

\begin{subequations}
	\begin{flalign}
	&\underset{\Theta^{R}}{\text{Min.}} \thinspace \sum_{a \in A} \sum_{i\in I}\left(C^{\text{R}+}_{i} R^{+}_{ia} + C^{\text{R}-}_i R^{-}_{ia}\right) \\
	& \text{subject to} \nonumber \\
	& \sum_{a} R_{ia}^{+} \le R^{+,\text{max}}_{i}, \quad \forall i \label{eq:A-Res-conv-max} \\
	& \sum_{a} R_{ia}^{-} \le R^{-,\text{max}}_{i}, \quad \forall i \label{eq:A-Res-conv-min} \\
	& \sum_{a} \left( R_{ia}^{+} + R_{ia}^{-} \right) \le P^{\text{max}}_{i}, \quad \forall i \label{eq:A-Res-conv-cap} \\
	& \sum_{i} R_{ia}^{+} \ge RR_a^{+}, \quad \forall a \label{eq:A-Res-req-up}\\
	& \sum_{i} R_{ia}^{-} \ge RR_a^{-}, \quad \forall a \label{eq:A-Res-req-dn} \\
	& \sum_{i \in \mathcal{M}_{b}^I} R_{ia}^{+} \le \sum_{l \in M_{a b}^{l}} X_l f^{\text{max}}_{l}, \quad \forall a, \forall b \ne a \label{eq:A-Res-Trans-up} \\
	& \sum_{i \in \mathcal{M}_{b}^I} R_{ia}^{-} \le \sum_{l \in M_{a {b}}^{l}} X_l f^{\text{max}}_{l}, \quad \forall a, \forall {b} \ne a \label{eq:A-Res-Trans-dn} \\
	& R^{+}_{ia},R^{-}_{ia}\geq 0, \forall i, \forall a
	\end{flalign}
\end{subequations}
The set of optimization variables $\Theta^{R}=\big\{{R}^{+}_{ia},{R}^{-}_{ia},  \enspace \forall i, a \big\}$.\\


\begin{figure}[h]
	\centering
	\includegraphics[scale=.6]{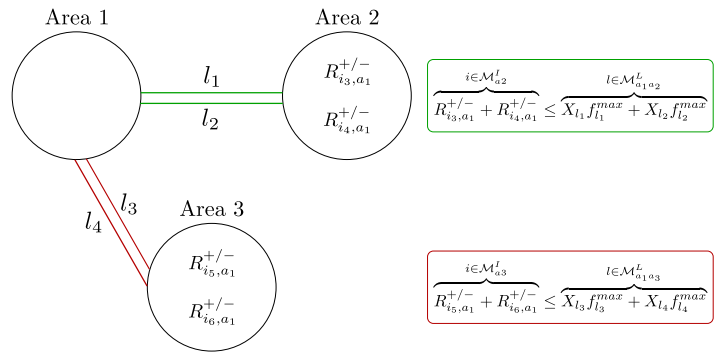}
	\caption{Illustrative example: Capacity allocation for reserve exchange for $a=\{a_1\}$, $b=\{a_2,a_3\}$ - constraints \eqref{eq:A-Res-Trans-up} and \eqref{eq:A-Res-Trans-dn}.}
	\label{fig:Wpen}
\end{figure}


\subsection*{Day-ahead market-clearing: Energy only}

\begin{subequations}
\begin{flalign}
&\underset{\Theta^{DA}}{\text{Min.}} \thinspace \sum_{i\in I}{C}_{i}{p}_{i} \\
& \text{subject to} \nonumber \\	
& \sum_{j \in \mathcal{M}_n^J} p_j^\text{W} + \sum_{i\in \mathcal{M}_n^I}{p}_{i} - D_n - \sum_{l \in L^{AC}\cup L^{DC}}  A_{ln} \left( \hat{f}^{\text{AC}}_{l} + \hat{f}_l^{DC} \right) = 0, \quad \forall n  \label{eq:A-CV-da-bal-constr} \\
& p_i \leq P^{\text{max}}_{i}-\sum_{a} R_{ia}^{+},  \quad \forall i  \label{eq:A-CV-da-conv-max}\\
& \sum_{a} R_{ia}^{-} \le p_i,  \quad \forall i  \label{eq:A-CV-da-conv-min}\\
& p_j^\text{W} \leq \bar{P}^{\text{W}}_{j},  \quad \forall j  \label{eq:A-CV-da-wind-max}\\
& \hat{f}^{\text{AC}}_{l} = B_l \sum_n A_{ln} \hat{\delta}_n, \quad \forall l \in L^{AC} \label{eq:A-CV-da-flowAC-constr}\\
& -(1-X_l) f^{\text{max}}_{l} \le \hat{f}^{\text{AC}}_{l} \le (1-X_l) f^{\text{max}}_{l}, \quad \forall l \in L^{AC} \label{eq:A-CV-da-flowAC-max}\\			
& -(1-X_l) f^{\text{max}}_{l} \le \hat{f}^{DC}_l \le (1-X_l) f^{\text{max}}_{l}, \quad \forall l \in L^{DC} \label{eq:A-CV-da-flowDC-constr}\\
& \hat{\delta}_1 = 0\\
& p_i \geq 0 , \forall i; \enspace p_j^\text{W} \geq 0 , \forall j; \nonumber \\ 
& \hat{\delta}_n \thinspace \text{free}, \forall n; \hat{f}_l^{AC} \thinspace \text{free}, \forall l \in  L^{AC}; \hat{f}_l^{DC} \thinspace \text{free}, \forall l \in  L^{DC}
\end{flalign}
\end{subequations}
The set of optimization variables $\Theta^{DA}=\big\{p_i,p_{j}^\text{W}, \hat{\delta}_n, \hat{f}^{AC}_l, \hat{f}^{DC}_l,  \enspace \forall i, j, n, l \big\}$.

\newpage

\subsection*{Balancing market-clearing}

\begin{subequations}
	\begin{flalign}
	&\underset{\Theta^{B}}{\text{Min.}} \thinspace \sum_{i \in I}{C}_{i}\left({r}_{i}^{+}-{r}_{i}^{-}\right)+\sum_{n \in N} {C}^{\text{sh}}l^{\text{sh}}_{n} \\
	& \text{subject to} \nonumber \\
	& \sum_{i \in \mathcal{M}_n^I} \left(r_{i}^{+}-r_{i}^{-} \right) + l^{\text{sh}}_n + \sum_{j \in \mathcal{M}_n^J}\left(P_{j}^*-p^\text{W}_j-p^{\text{spill}}_{j}\right) \nonumber \\
	& + \sum_{l \in L^{AC}\cup L^{DC}}  A_{ln} \left(\hat{f}^{\text{AC}}_{l} -\tilde{f}^{\text{AC}}_{l} + \hat{f}^{DC}_{l} - \tilde{f}^{DC}_{l} \right) = 0, \quad  \forall n   \label{eq:A-CV-rt-bal-constr} \\	
	& \tilde{f}^{\text{AC}}_{l} = B_l \sum_n A_{ln} \tilde{\delta}_{n}, \quad \forall l \in L^{AC}  \label{eq:A-CV-rt-flowAC-constr}\\
	& r^+_{i} \leq \sum_{a} R_{ia}^{+}, \quad \forall i \label{eq:A-CV-rt-conv-res-up-max}\\
	& r^-_{i} \leq \sum_{a} R_{ia}^{-}, \quad \forall i \label{eq:A-CV-rt-conv-res-dn-max}\\
	& - f^{\text{max}}_{l} \le \tilde{f}^{\text{AC}}_{l} \le f^{\text{max}}_{l}, \quad \forall l \in L^{AC}  \label{eq:A-CV-rt-flowAC-max}\\
	& -f^{\text{max}}_{l} \le \tilde{f}^{DC}_{l} \le f^{\text{max}}_{l}, \quad \forall l \in L^{DC}  \label{eq:A-CV-rt-flowDC-constr}\\
	& l^{\text{sh}}_{n} \leq D_n, \quad \forall n   \label{eq:A-CV-rt-load-shed-max}\\
	& p^{\text{spill}}_{j} \leq P^*_{j}, \quad \forall j  \label{eq:A-CV-rt-wind-spill-max}\\	
	& \tilde{\delta}_{1} = 0 \label{eq:A-CV-ref-node}\\
	&  r^+_{i},r^-_{i}\geq 0, \forall i;\enspace p^{\text{spill}}_{j} \geq 0, \forall j;  \enspace l^{\text{sh}}_{n} \geq 0, \forall n; \nonumber \\ 
	& \tilde{\delta}_{n} \thinspace \text{free}, \forall n; \tilde{f}_{l}^{AC} \thinspace \text{free}, \forall l \in  L^{AC}; \enspace \tilde{f}_{l}^{DC} \thinspace \text{free}, \forall l \in  L^{DC}		
	\end{flalign}
\end{subequations}
The set of optimization variables $\Theta^{B}=\big\{{r}^{+}_{i},{r}^{-}_{i}, p^{\text{spill}}_{j}, {l}^{\text{sh}}_{n}, \tilde{\delta}_{n}, \tilde{f}^{AC}_{l}, \tilde{f}^{DC}_{l}, \enspace \forall i, j, n, l  \big\}$.\\


Constraints (13c) and (13d) imply that reserve deployment depends only on reserve availability in real-time, irrespective of the initial area requirements. There is a single pool of balancing resources accessible from all areas during real-time operation.

}

\end{document}